\documentclass[11pt]{amsart}
\usepackage{amssymb}
\usepackage{caption}
\usepackage{subcaption}
\usepackage{graphicx}
\usepackage{xcolor} 
\usepackage{tensor}
\usepackage{fullpage} 
\usepackage{amsmath}
\usepackage{amsthm}
\usepackage{verbatim}
\usepackage{hyperref}
\usepackage{enumitem}
\usepackage{mathtools}
\providecommand{\bysame}{\leavevmode\hbox to3em{\hrulefill}\thinspace}
\providecommand{\MR}{\relax\ifhmode\unskip\space\fi MR }
\providecommand{\MRhref}[2]{%
	\href{http://www.ams.org/mathscinet-getitem?mr=#1}{#2}
}
\providecommand{\href}[2]{#2}
\setlist[enumerate]{leftmargin=1.5em}
\setlist[itemize]{leftmargin=1.5em}

\providecommand{\bysame}{\leavevmode\hbox to3em{\hrulefill}\thinspace}
\providecommand{\MR}{\relax\ifhmode\unskip\space\fi MR }
\providecommand{\MRhref}[2]{%
	\href{https://urldefense.com/v3/__http://www.ams.org/mathscinet-getitem?mr=*1*7D*7B*2__;IyUlIw!!Mih3wA!Ab5Z3hIsC2h98i8YRQ7VOKwfACgAgCvQF4vgVmE6EG5Lyy5uaKCmO6MHj8j0Xuk0aNaX_fziDACCPD0EYw$  }
}
\providecommand{\href}[2]{#2}

\usepackage{seqsplit} 
\definecolor{green}{rgb}{0,0.8,0} 

\newtheorem{theorem}{Theorem}[section]
\newtheorem{corollary}[theorem]{Corollary}

\newtheorem{proposition}[theorem]{Proposition}

\theoremstyle{definition}

\theoremstyle{remark}
\newtheorem{remark}[theorem]{Remark}

\numberwithin{equation}{section}

\newcommand{\abs}[1]{\left\vert#1\right\vert}

\newcommand{\set}[1]{\left\{#1\right\}}

\newcommand{\nnrm}[1]{{\vert\kern-0.25ex\vert\kern-0.25ex\vert #1 
		\vert\kern-0.25ex\vert\kern-0.25ex\vert}}

\newcommand{\lap}{\Delta}

\newcommand{\rd}{\partial}
\newcommand{\nb}{\nabla}

\newcommand{\alp}{\alpha}
\newcommand{\bt}{\beta}
\newcommand{\gmm}{\gamma}

\newcommand{\dlt}{\delta}

\newcommand{\eps}{\epsilon}

\newcommand{\kpp}{\kappa}

\newcommand{\tht}{\theta}

\newcommand{\Omg}{\Omega}

\newcommand{\bfN}{{\bf N}}

\newcommand{\bfT}{{\bf T}}

\newcommand{\bbR}{\mathbb R}

\newcommand{\bbT}{\mathbb T}

\begin{document}
	
	\title{On evolution of corner-like gSQG patches}
	\author{Junekey Jeon}
	\address{Department of Mathematics, UC San Diego.}
	\email{j6jeon@ucsd.edu}
	\author{In-Jee Jeong}
	\address{Department of Mathematics and RIM, Seoul National University.}
	\email{injee\_j@snu.ac.kr}
	\date{\today}
	
	
	
	\begin{abstract}
		We study the evolution of corner-like patch solutions to the generalized SQG equations. Depending on the angle size and order of the velocity kernel, the corner instantaneously bents either downward or upward. In particular, we obtain the existence of strictly convex and smooth patch solutions which become immediately non-convex. 
	\end{abstract}
	
	\maketitle
	

	
	\section{Introduction}
	
	\subsection{Patch solutions}
	
	We consider patch solutions to the generalized surface quasi-geostrophic (gSQG) equations, which are given by \begin{equation}\label{eq:gSQG}
		\left\{
		\begin{aligned}
			\rd_t \tht + u \cdot \nb \tht = 0, & \\
			u = -\nb^\perp (-\lap)^{-1+\frac{\alp}{2}} \tht ,
		\end{aligned}
		\right.
	\end{equation} where $\tht(\,\cdot\,,t):\bbR^2\to\bbR$ and $u(\,\cdot\,,t):\bbR^2\to\bbR^2$. Here, $0\le\alp<2$ is a parameter, and the special cases $\alp=0, 1$ correspond to the two-dimensional incompressible Euler and SQG equations, respectively. In this paper, by \emph{a patch solution} to \eqref{eq:gSQG} we mean the indicator function
	$\mathbf{1}_{\Omega(t)}$ over a moving bounded open set $\Omega(t)\subseteq\mathbb{R}^{2}$,
	whose boundary is given by a simple closed curve $\gamma(\,\cdot\,,t)\colon\mathbb{T}\to\mathbb{R}^{2}$ satisfying the following \emph{contour dynamics equation}:
	\begin{equation}\label{eq:CDE}
		\partial_{t}\gamma(\xi,t) =
		\begin{dcases}
			c_{\alpha}\int_{\mathbb{T}}\partial_{\xi}\gamma(\eta,t)
			\ln\abs{\gamma(\xi,t) - \gamma(\eta,t)}\,d\eta
			+ \lambda(\xi,t)\partial_{\xi}\gamma(\xi,t)
			&\textrm{if $\alpha = 0$,} \\
			c_{\alpha}\int_{\mathbb{T}}
			\frac{\partial_{\xi}\gamma(\xi,t) - \partial_{\xi}\gamma(\eta,t)}
			{\abs{\gamma(\xi,t) - \gamma(\eta,t)}^{\alpha}}\,d\eta
			+ \lambda(\xi,t)\partial_{\xi}\gamma(\xi,t)
			&\textrm{otherwise,}
		\end{dcases}
	\end{equation}
	where $c_{\alpha}>0$ is a constant and $\lambda$ is any scalar function that can depend on $\gamma$. It is shown in \cite{Gan,KYZ} that as long as the curve $\gamma$ solves \eqref{eq:CDE}, remains to be a simple closed curve, and retains certain regularity, then the function $(x,t)\mapsto\mathbf{1}_{\Omega(t)}(x)$ defines a weak solution to \eqref{eq:gSQG}. (The mentioned references only talk about the case $0<\alpha\leq 1$ but the argument extends to the whole range as well.) Moreover, with a careful choice of the function $\lambda$, it is known that \eqref{eq:CDE} admits a unique solution locally in time if the initial patch boundary is sufficiently regular.
	
	In the specific case of $\alp=0$ (the Euler equations), taking $\lambda\equiv 0$ gives that \eqref{eq:CDE} is nothing but the characteristic equation restricted to the patch boundary, which is known to have a unique solution by Yudovich theory (\cite{Y1}). The resulting patch solution is also known to exists globally in time and preserves its boundary regularity (\cite{Chemin1993,BeCo,Ser1}). When $0<\alpha<1$, the characteristic equation does not necessarily admit a unique solution, but \eqref{eq:CDE} still does with $\lambda\equiv 0$ (\cite{Gan}). When $1\leq \alpha<2$, the velocity field $u$ may not be even defined on the boundary, but if we select
	\begin{equation}\label{eq:lambda}\begin{split}
			\lambda(\xi,t) &= c_{\alpha}\frac{\xi+\pi}{2\pi}
			\int_{\mathbb{T}}
			\frac{\partial_{\xi}\gamma(\eta,t)}
			{\abs{\partial_{\xi}\gamma(\eta,t)}^{2}}\cdot
			\partial_{\eta}\left(
			\int_{\mathbb{T}}\frac{\partial_{\xi}\gamma(\eta,t) - \partial_{\xi}\gamma(\zeta,t)}
			{\abs{\gamma(\eta,t) - \gamma(\zeta,t)}^{\alpha}}\,d\zeta
			\right)d\eta \\
			&\quad\quad\quad\quad
			- c_{\alpha}\int_{-\pi}^{\xi}
			\frac{\partial_{\xi}\gamma(\eta,t)}
			{\abs{\partial_{\xi}\gamma(\eta,t)}^{2}}\cdot
			\partial_{\eta}\left(
			\int_{\mathbb{T}}\frac{\partial_{\xi}\gamma(\eta,t) - \partial_{\xi}\gamma(\zeta,t)}
			{\abs{\gamma(\eta,t) - \gamma(\zeta,t)}^{\alpha}}\,d\zeta
			\right)d\eta,
	\end{split}\end{equation}
	then \eqref{eq:CDE} still admits a unique solution (\cite{GaPa,CCCGW,CCG-TRAN,Gan}). From now on, we implicitly assume that $\lambda$ is chosen to be identically zero when $0\leq\alpha<1$, or as in \eqref{eq:lambda} if $1\leq\alpha<2$.
	
	\subsection{Main results}
	
\begin{figure}
		\centering
		\includegraphics[scale=0.5]{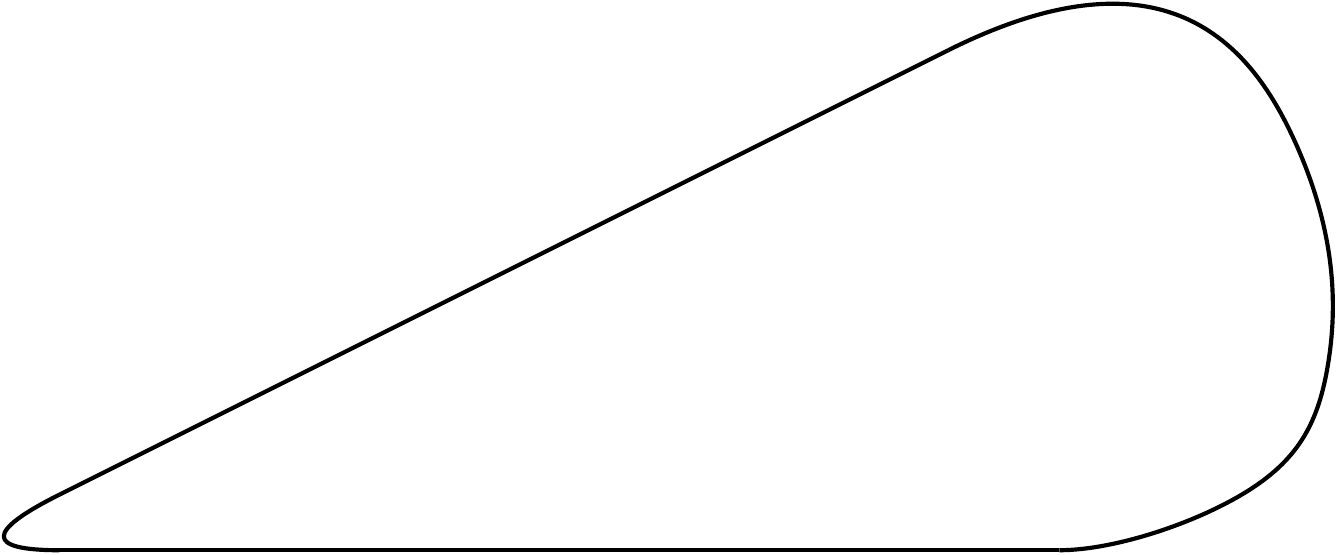}  
		\caption{An example of domain $\Omg_0$.} \label{fig:corner}
	\end{figure}

	In the current work, we are interested in the evolution of patch solutions which locally look like a \textit{corner} at the initial time. To be more precise, we assume that the initial patch $\Omg_0$ has two length scales $0<\dlt\ll M$ such that \begin{equation}\label{eq:corner}
		\begin{split}
			\Omg_0 \cap \left\{ x \in [-M,M]^{2} \colon \abs{x} > \delta \right\}
			= \left\{ x \in (0,\infty)^{2} \colon \beta x_2 < x_1 \right\}
			\cap \left\{ x \in [-M,M]^{2} \colon \abs{x} > \delta \right\}
		\end{split}
	\end{equation} for some $0\leq \bt < \infty$, see Figure~\ref{fig:corner}. For simplicity, we assume the following on the area $|\Omg_0|$: \begin{equation}\label{eq:corner2}
		\begin{split}
			|\Omg_0| \le 100M^2. 
		\end{split}
	\end{equation} The case $\bt=0$ corresponds to corners with the right angle. We call the intersections of $\partial\Omega_{0}$ with the line segments $\set{x\in(0,M)^{2}\colon \abs{x}>\delta,\ \beta x_{2}=x_{1}}$ and $(\delta,M)\times\set{0}$ as \emph{the upper edge} and \emph{the lower edge} of $\Omega_{0}$, respectively.

	It turns out that generically, during the evolution, parts of both of the edges of $\Omega_{0}$ sufficiently close to the origin instantaneously ``bents'' either downward or upward. Our main result gives a precise criterion on $\alp$ and $\bt$ that determines which one of the two is the case.
	\begin{theorem}\label{thm:main}
		Assume that we are given a patch $\Omg_0$ with smooth boundary satisfying \eqref{eq:corner}--\eqref{eq:corner2}. 
		
		Then, the boundary of the corresponding unique patch solution $\Omg(t)$ becomes instantaneously convex downward (upward, resp.) in the region $A \cap \Omg(t)$ if \begin{equation}\label{eq:convex-criterion}
			\begin{split}
				F(\alp,\bt) := \frac{\beta}{\beta^{2}+1}
				- \frac{\alpha}{(\beta^{2}+1)^{1-\alpha/2}}
				\int_{-\tan^{-1}\beta}^{\pi/2}\cos^{\alpha}\theta\,d\theta > 0 \qquad (\mbox{resp.} < 0), 
			\end{split}
		\end{equation} where the annular region $A$ is defined as
		\begin{equation}\label{eq:annulus}
			\begin{split}
				A\coloneqq \left\{ x \in \bbR^2 : {C_\alp^{-1}}|F(\alp,\bt)|^{-\frac12} \dlt < |x| < C_\alp |F(\alp,\bt)|^{\frac{1}{1+\alp}}M \right\}
			\end{split}
		\end{equation}
		with some $C_\alp>0$ depending only on $\alp$. 
	\end{theorem}
	
	\begin{figure}
		\centering
		\includegraphics[scale=0.8]{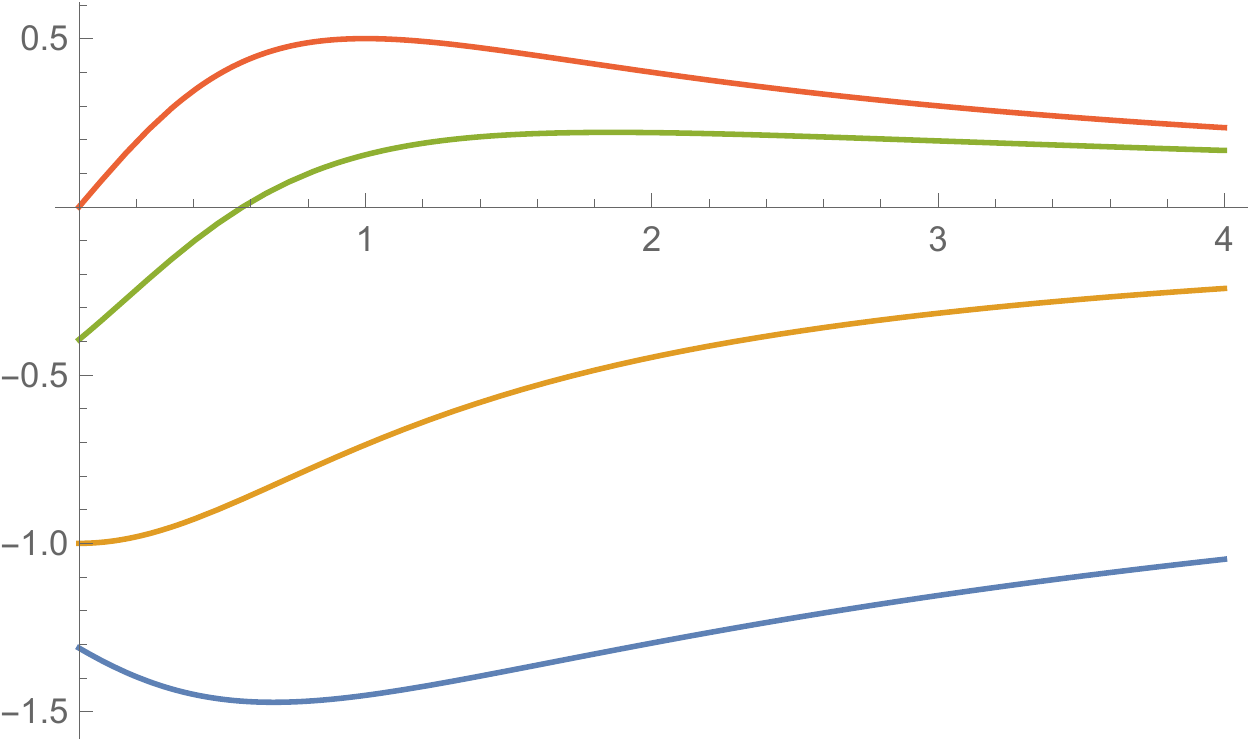}  
		\caption{Plot of $F(\alp,\cdot)$ in the region $0\le\bt\le4$ with $\alp = 0,  0.5, 1, 1.5$ (from top to bottom).} \label{fig:plot}
	\end{figure}
	
\begin{figure}
		\centering
		\begin{subfigure}{0.48\textwidth}
			\includegraphics[width=\textwidth, trim=180 0 180 0]{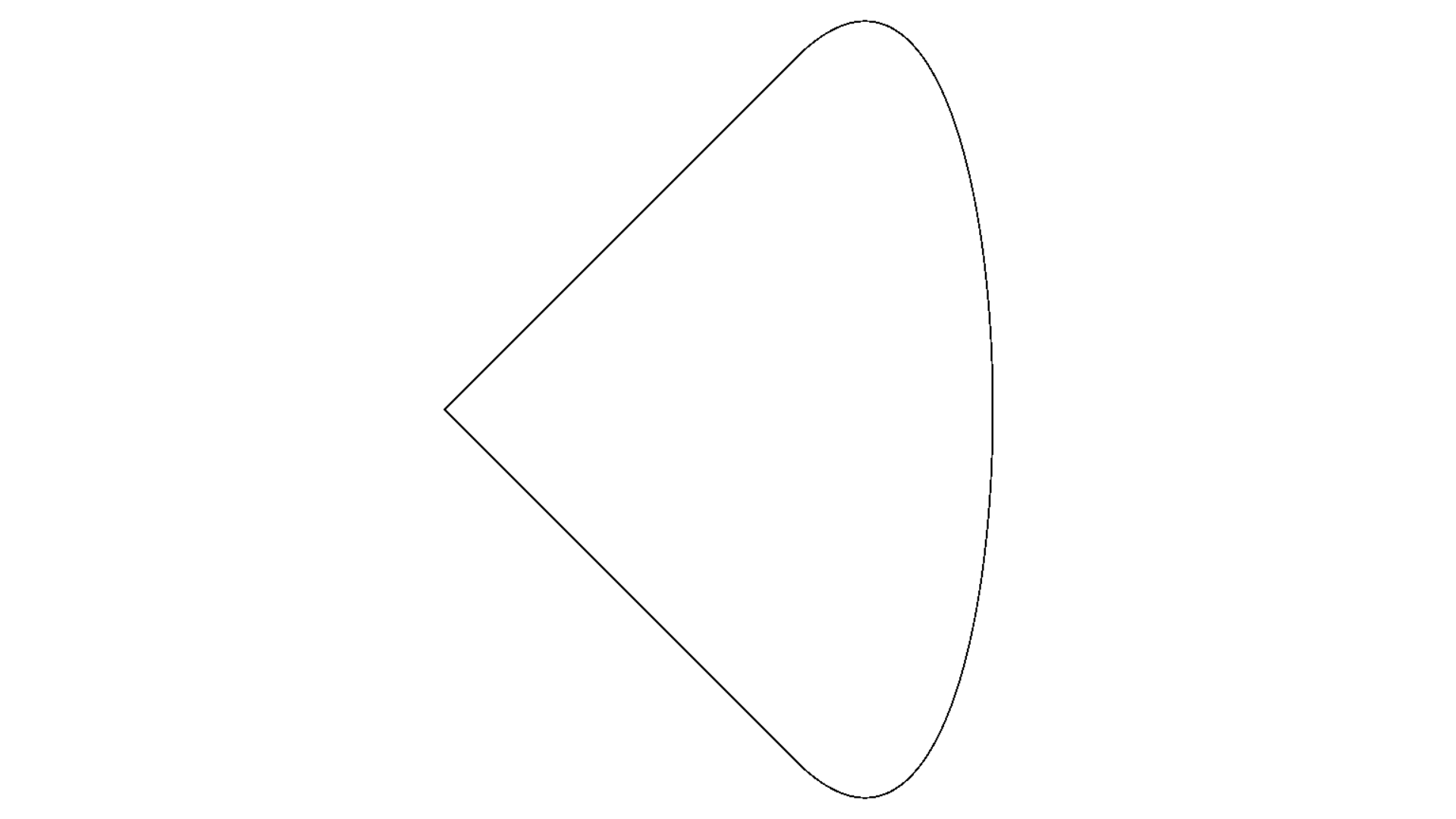}
			\includegraphics[width=\textwidth, trim=180 0 180 0]{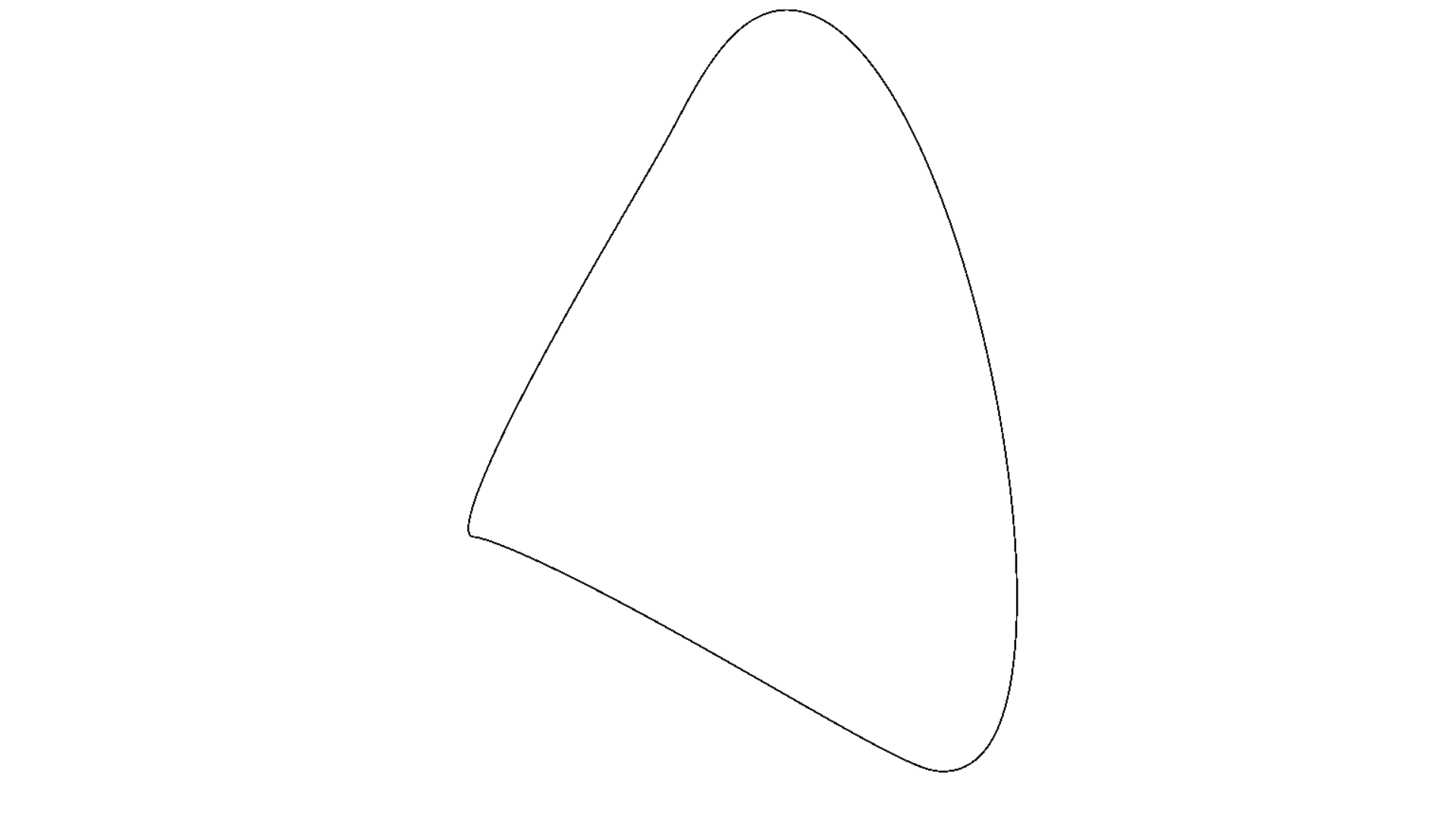}
			\includegraphics[width=\textwidth, trim=180 0 180 0]{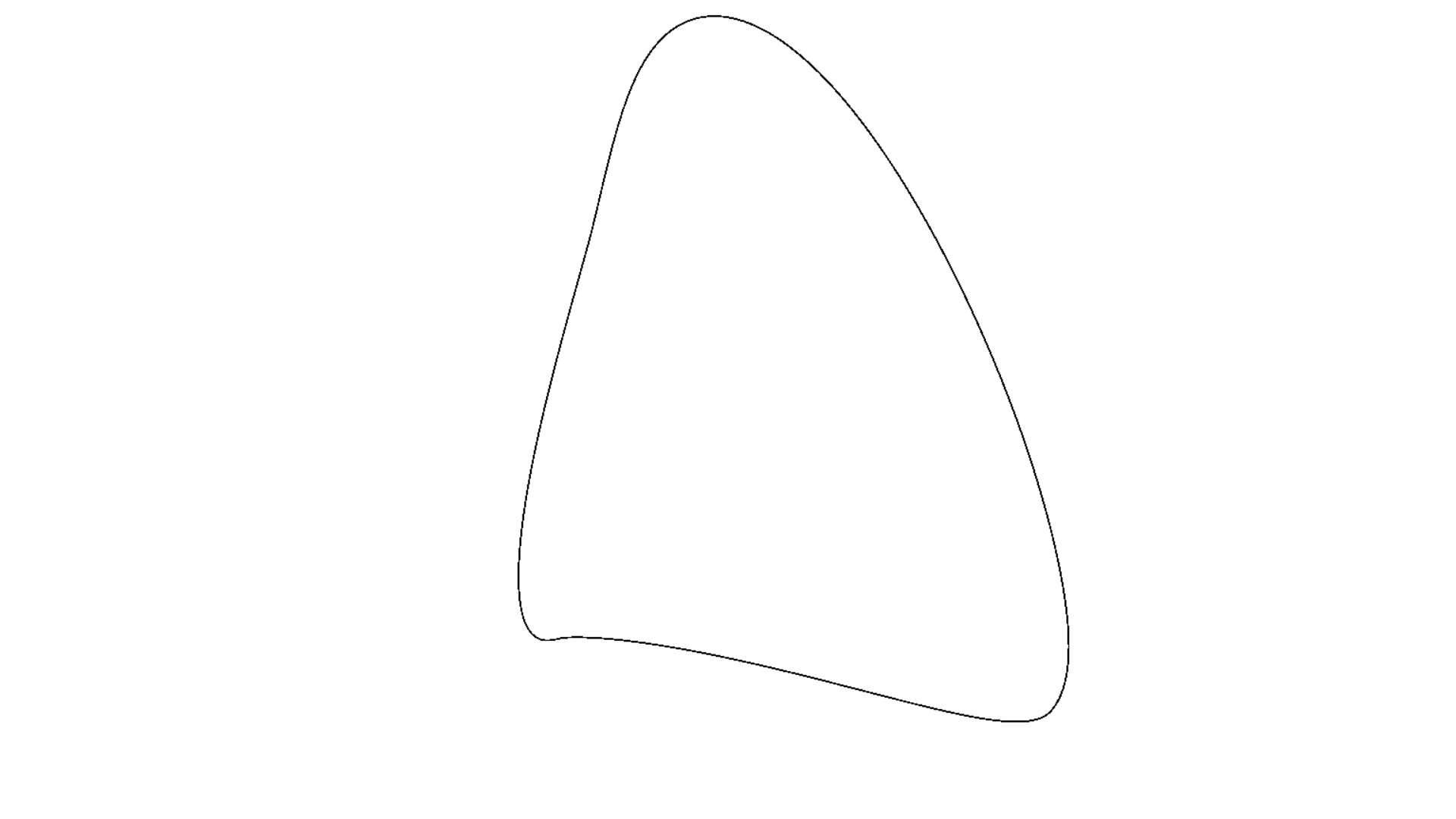}
			\caption{$\beta=0$}
		\end{subfigure}
		\hfill
		\begin{subfigure}{0.48\textwidth}
			\includegraphics[width=\textwidth, trim=180 0 180 0]{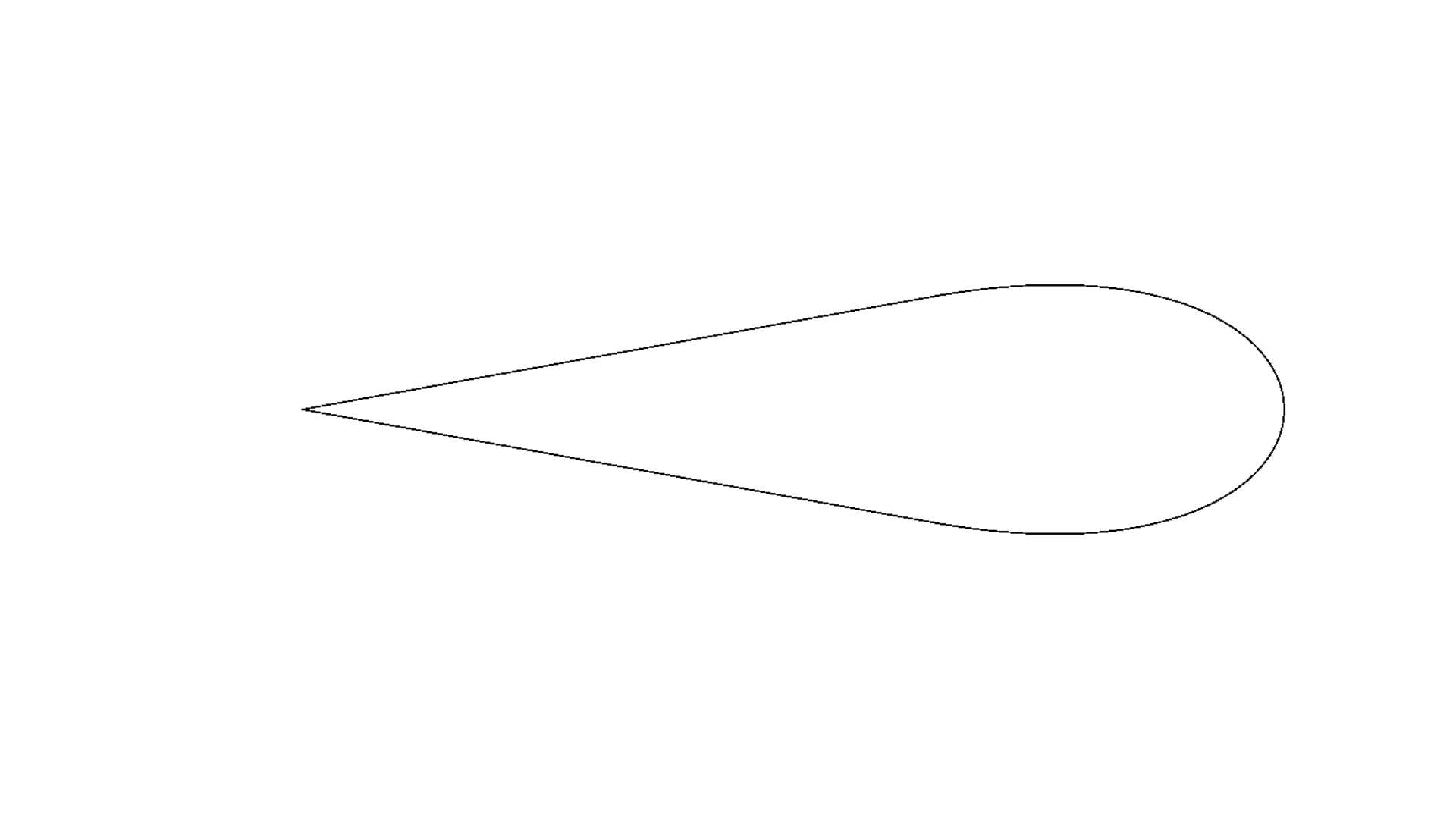}
			\includegraphics[width=\textwidth, trim=180 0 180 0]{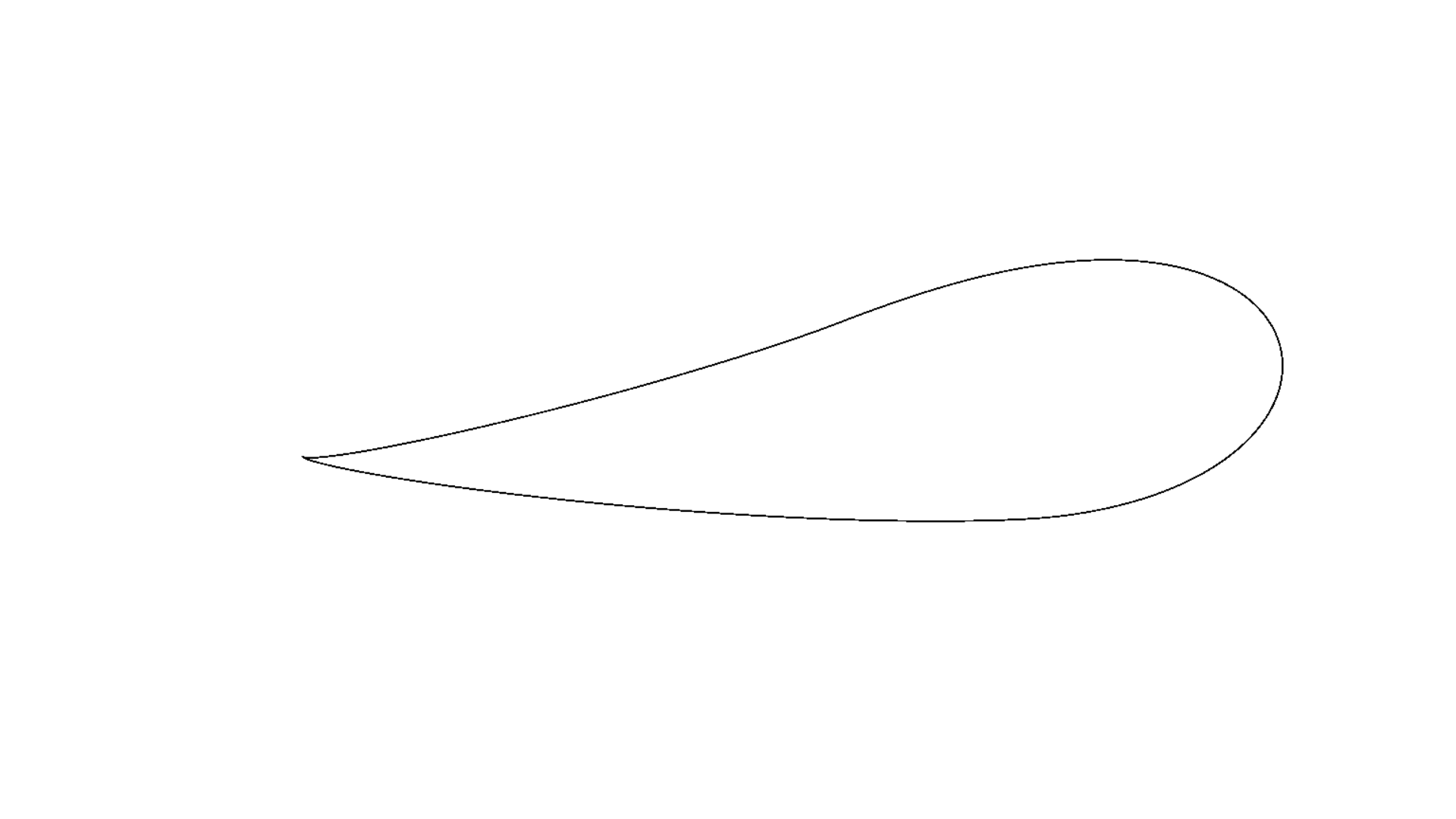}
			\includegraphics[width=\textwidth, trim=180 0 180 0]{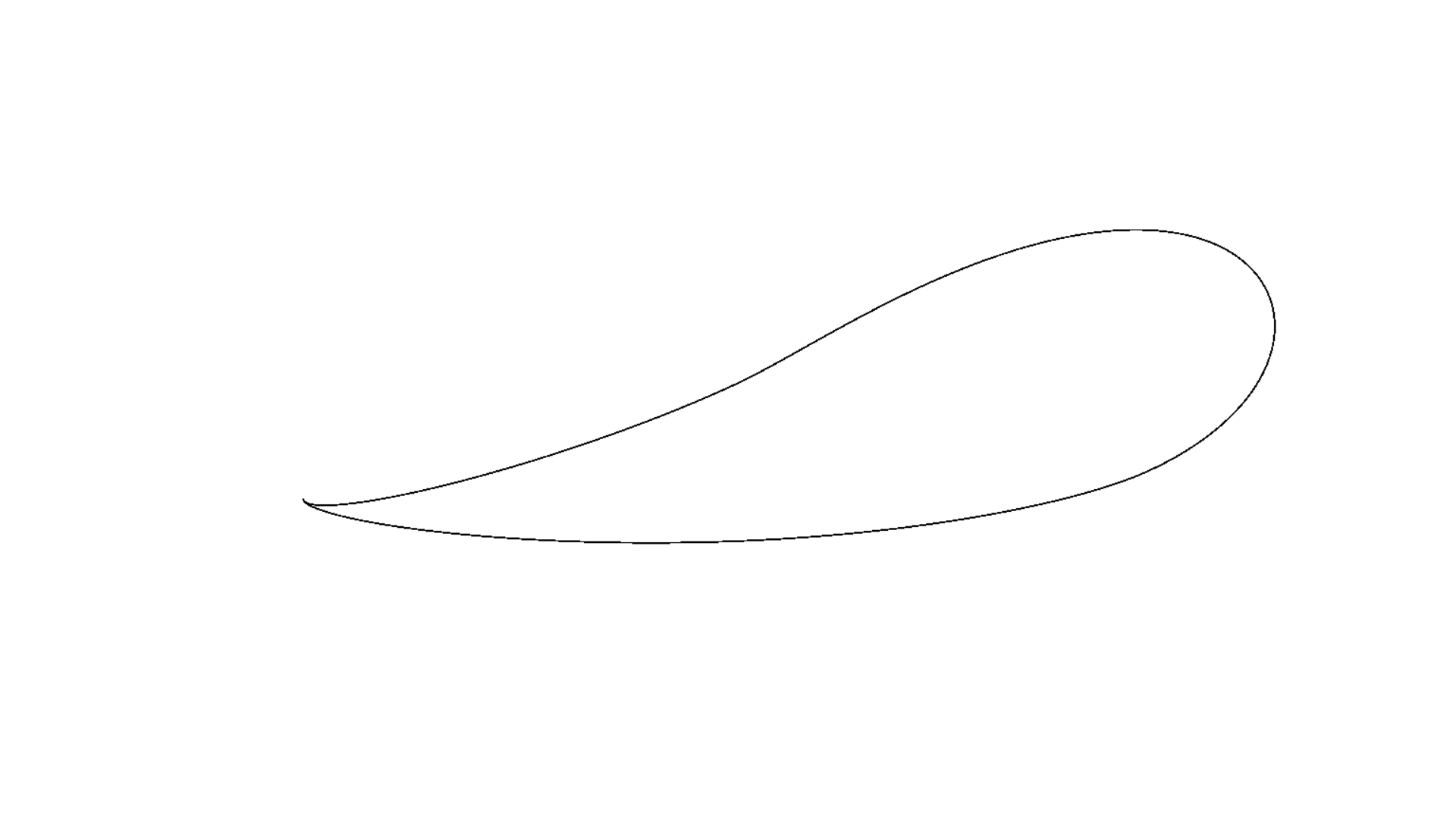}
			\caption{$\beta=2.7$}
		\end{subfigure}
		\caption{Evolution of corner-like patches for $\alpha=0.3$ at $t=0, 0.2$, and $0.4$.}
		\label{fig:simul}
	\end{figure}

	\begin{remark}
		As one can see from Figure~\ref{fig:plot}, $F(\alp,\bt)$ takes both positive and negative signs. Indeed, one can compute that
		\[
		\frac{\partial}{\partial\beta}
		\left((\beta^{2}+1)^{1-\alpha/2}F(\alpha,\beta)\right)
		= -\frac{\alpha-1}{(\beta^{2}+1)^{\alpha/2}}
		\] holds. Based on this formula, we have the following:
		\begin{itemize}
			\item In the Euler case $\alp=0$, we have $F(0,\bt)>0$ for all $\bt>0$, which means that for any angle of the corner in the range $\left(0,\frac{\pi}{2}\right)$, the patch boundary near the corner always bents downward. This is consistent with the numerical simulations from \cite{CS} (their sign convention for the velocity is the opposite from ours).
			
			\item In the SQG case $\alp=1$ and the more singular case $\alp>1$, we have $F(\alp,\bt) < 0$ for all $\bt\geq 0$, which means that for any angle of the corner in the range $\left(0,\frac{\pi}{2}\right]$, the patch boundary near the corner always bents upward.
			
			\item When $0<\alpha<1$ (corresponding to the equations interpolating Euler and SQG), $F(\alpha,0)<0$ and $\lim_{\beta\to\infty}(\beta^{2}+1)^{1-\alpha/2}F(\alpha,\beta) = \infty$, so $(\beta^{2}+1)^{1-\alpha/2}F(\alpha,\beta)$ strictly monotonically increases from some negative value to $\infty$. In particular, there is a unique zero $\bt^*(\alp)$ of $F(\alp,\,\cdot\,)$. Consequently, the patch boundary near the corner bents upward for large corner angle ($\beta<\beta^{*}(\alpha)$) while it bents downward for small corner angle ($\beta>\beta^{*}(\alpha)$). This is consistent with our own numerical simulations shown in Figure~\ref{fig:simul}.
		\end{itemize} 
	\end{remark}
	
	As a consequence of the above, we obtain the existence of initially convex and smooth patches which become non-convex at a later time. To the best of our knowledge, such a result has not been proved even in the case of Euler. 
	\begin{corollary}\label{cor:main}
		For any $0\le \alp< 2$, there exists a strictly convex patch $\Omg_0$ with smooth boundary such that the corresponding unique patch solution becomes non-convex immediately after the initial time.
	\end{corollary}

	\subsection{Background}
	
	To put the above results into context, let us review relevant works which have inspired our investigation. To begin with, it is a natural question to ask what happens to patches with a corner for the gSQG equations, since in contrast to the local well-posedness results in the smooth patch regime, not much is known without the smoothness assumption. 
	
	In the exceptional case of the Euler equations where $\alp=0$, there is an existence and uniqueness result for patch solutions without a smoothness assumption (\cite{Y1}). When the initial patch has a corner, numerical computations of \cite{CS} and formal asymptotic analysis of \cite{EJSVP1} suggest that it instantaneously evolves into a cusp, although rigorously proving it seems to be an open problem. Our main theorem can be applied \textit{without} the local smoothing parameter $\dlt$ in the Euler case, and it shows that the corner patch immediately bents downward for any corner angle strictly between $0$ and $\pi/2$, confirming the results of \cite{CS}. In particular, we deduce that the patch is not convex for any small enough $t$.
	
	The main step in our proof is to compute the precise asymptotic behavior of the normal velocity field associated with the corner-type patch. Then, the change in convexity of the patch boundary can be derived from the curvature evolution equation in the very recent work of Kiselev--Luo \cite{KiLuo2}. A particularly interesting feature revealed by our computation is that in the case of $0\leq\alp<1$, there is a special value of $\bt^* = \bt^*(\alp)$ such that the corner evolution is the ``opposite'' for $\bt<\bt^*$ and $\beta>\bt^*$. In fact, it turns out that in the regime $M\to\infty$ and $\delta\to 0$ (that is, when the corner is infinitely long and perfectly sharp), the normal velocity on both of the edges becomes a constant if the corner angle is equal to $\beta^{*}$; see Remark~\ref{remk:velocity estimate}. This shows that, formally, infinitely long corners with angle $\bt^*$ are steady solutions for the gSQG equations relative to the frame moving together with the corner. In particular, we have $\bt^*(0)=0$, which is consistent with the (conjectured) existence of relative equilibria of Euler patches containing $90^{\circ}$ corners (\cite{OverII,HMW,HHHM,ZOWZ,WOZ}). 
	
	
	In general, numerical simulations on patch solutions have revealed very complicated dynamics of the boundary, involving repeated filamentation and fast curvature growth (\cite{CS,Dr,DrMc,Scott,SD}). Most notably, a recent work of Scott and Dritschel (\cite{SD}) shows the possibility of locally self-similar finite time singularity formation for gSQG patches. Very recently, the existence of self-similar solutions which exhibit similar behavior was proved by Garc\'{i}a--G\'{o}mez-Serrano \cite{GSerrano2}. When the domain is the upper half-plane, rapid small scale creation for the patch boundary including finite time blow-up for patch solutions have been established (\cite{KiLi,KYZ,KRYZ,GaPa}), thanks to the stability of the instability provided by the boundary. To the best of our knowledge, such a dramatic growth for sufficiently regular patches has not been proved to exist for domains without boundary like $\bbR^2$ or $\bbT^2$, but see \cite{CJ-AML,CJ-NA,EJSVP2} for some results in this direction. On the other hand, in the aforementioned work \cite{KiLuo2}, the authors observe a remarkable structure in the curvature equation, and in particular applied it to prove strong ill-posedness of the patch problem with the boundary regularity of $C^2$ for 2D Euler. Namely, it is possible for a patch with a $C^{2}$ boundary (which in particular implies finite curvature) to instantaneously have infinite curvature for $t\ne0$.

	\section{Proofs}
	
	\subsection{Curvature evolution equation}
	
	We denote the unit tangent and normal vectors as $\bfT$ and $\bfN = -\bfT^\perp$, respectively, defined on the patch boundary. Then, we define the (signed) curvature by \begin{equation}\label{eq:curvature-def}
		\begin{split}
			\kpp = -\rd_s \bfT\cdot\bfN = \bfT\cdot\rd_s\bfN,
		\end{split}
	\end{equation} where $\rd_s$ is the differentiation with respect to the arc-length parameter $s$. Given a simple closed smooth curve $\gmm$ in $\bbR^2$, it is well-known that the domain enclosed by $\gmm$ is convex if and only if either $\kappa$ is non-positive or non-negative everywhere on the curve, depending on the orientation of the parametrization. We will always assume counterclockwise orientation, and in that case the domain is convex if and only if $\kappa\geq 0$ everywhere.
	
	Now, consider a parametrized simple closed curve $\gamma(\,\cdot\,,t)\colon\mathbb{T}\to\mathbb{R}^{2}$ evolving over time, whose evolution is governed by a velocity field $v(\,\cdot\,,t)\colon\mathbb{T}\to\mathbb{R}^{2}$ along the curve in the sense that
	\begin{equation}\label{eq:gmm-evol}
		\begin{split}
			\partial_{t}\gamma(\xi,t) = v(\xi,t)
		\end{split}
	\end{equation} 
	holds. Assuming that $v$ is at least twice continuously differentiable with respect to the curve parameter $\xi$, Kiselev and Luo \cite{KiLuo} have derived the evolution equation \begin{equation}\label{eq:kpp-evol}
		\begin{split}
			\rd_t\kpp = - 2\kpp \rd_s v\cdot\bfT- \rd_s^2 v\cdot \bfN
		\end{split}
	\end{equation} for the curvature $\kappa$ of the curve $\gamma$ evaluated at a fixed $\xi\in\mathbb{T}$, where $\rd_s$ can be computed using the formula  $\partial_{s}=\frac{1}{\abs{\partial_{\xi}\gamma(\xi,t)}}\partial_{\xi}$. Substituting the right-hand side of \eqref{eq:CDE} to $v$, we obtain the evolution equation for the curvature in the contour dynamics equation.
	
	\subsection{The main estimate and the proof of Theorem~\ref{thm:main}}

	We now state our main estimate, which implies Theorem~\ref{thm:main}.
	\begin{proposition}\label{prop:no strict convexity}
		Consider the gSQG patch solution $\Omega(\cdot)$ with the initial data $\Omega_0$ having smooth boundary and satisfying \eqref{eq:corner}--\eqref{eq:corner2}. Then, the corresponding velocity at $t=0$ satisfies the estimate \begin{equation}\label{eq:vel-second-grad}
			\begin{split}
				\abs{\partial_{s}^{2}v_{2}((a,0),0) - \frac{c_{\alpha}F(\alpha,\beta)}{a^{1+\alpha}}}
				\leq C_\alp \left( \frac{1}{M^{1+\alp}} + \frac{\dlt^{2}}{a^{3+\alpha}} \right)
			\end{split}
		\end{equation} for $2\dlt \le a \le M/2$, where $C_\alp>0$ is a constant depending only on $\alp$.
	\end{proposition}
	
	\begin{proof}
		Assume that we have chosen a parametrization $\gamma(\,\cdot\,,0)$ of $\partial\Omega(0) $ such that $\gamma(\pi,0) = (a,0)$ for some $2\delta\leq a\leq M/2$. We will derive a sharp condition on $\alp$ and $\bt$ which determines the sign of the right-hand side of \eqref{eq:kpp-evol} evaluated at $((a,0),0)$, when $\delta\ll a\ll M$.
		
		Note that the {vertical velocity} at $((a,0),0)$ is given as
		\[
		v_{2}(\pi,0) =
		c_{\alpha}\int_{\mathbb{T}}
		\frac{(\partial_{\xi}\gamma(\pi,0) - \partial_{\xi}\gamma(\eta,0))\cdot(0,1)}
		{\abs{(a,0) - \gamma(\eta,0)}^{\alpha}}\,d\eta
		= c_{\alpha}\mathrm{p.v.}\int_{\Omega_{0}^{\delta}}
		\frac{a - y_{1}}{\left((a-y_{1})^{2} + y_{2}^{2}\right)^{1 + \alpha/2}}
		\,dy_{1}\,dy_{2}.
		\]
		We shall obtain an asymptotic formula for the second derivative of the above with respect to $a$. Define
		\begin{align*}
			A(a) &\coloneqq \mathrm{p.v.}\int_{\Omega_{0}\cap[-M,M]^{2}}
			\frac{a - y_{1}}{\left((a-y_{1})^{2} + y_{2}^{2}\right)^{1 + \alpha/2}}
			\,dy_{1}\,dy_{2}, \\
			B(a) &\coloneqq \int_{\Omega_{0}\setminus[-M,M]^{2}}
			\frac{a - y_{1}}{\left((a-y_{1})^{2} + y_{2}^{2}\right)^{1 + \alpha/2}}
			\,dy_{1}\,dy_{2}, \quad\textrm{and} \\
			C(a) &\coloneqq \int_{\Omega_{0}\setminus\Omega_{0}^{\delta}}
			\frac{y_{1} - a}{\left((a-y_{1})^{2} + y_{2}^{2}\right)^{1 + \alpha/2}}
			\,dy_{1}\,dy_{2},
		\end{align*}
		so that $v_{2}(\pi,0) = c_{\alpha}(A(a) + B(a) + C(a))$.
		
		We estimate $B''(a)$ first. From the expression
		\begin{equation*}
			B'(a) = \int_{\Omg_0  \setminus [-M,M]^2 }
			\frac{1}{\left((a-y_{1})^{2}+y_{2}^{2}\right)^{1+\alpha/2}}
			- \frac{(\alpha+2)(a-y_1)^2}{\left((a-y_{1})^{2}+y_{2}^{2}\right)^{2+\alpha/2}}
			\,dy_1 \, dy_2,
		\end{equation*}
		we easily derive
		\begin{align*}
			|B'(a)| & \le \int_{\Omg_0  \setminus [-M,M]^2 } \frac{C_{\alpha}}{(a-y_1)^{2+\alpha} }
			\,dy_1 \, dy_2
			\leq \frac{C_{\alpha}}{(M-M/2)^{2+\alpha}}\abs{\Omega_{0}}
			\leq \frac{C_{\alpha}}{M^{2+\alpha}}\abs{\Omega_{0}}
		\end{align*}
		for some constant $C_{\alpha}$ that may change from side to side. In a similar manner, one can show
		\begin{equation}\label{eq:B estimate}
			|B''(a)| \leq \frac{C_{\alpha}}{M^{3+\alpha}}\abs{\Omg_0} .
		\end{equation} 
		
		In the same way, we obtain
		\begin{equation}\label{eq:C estimate}
			|C''(a)| \leq \int_{B_{0}(\delta)}\frac{C_{\alpha}}{\abs{a-y_{1}}^{3+\alpha}}\,dy_{1}\,dy_{2}
			\leq \frac{C_{\alpha}}{a^{3+\alpha}}\delta^{2}.
		\end{equation}
		
		Next, we rewrite $A(a)$ as
		\begin{align*}
			A(a) &= \int_{0}^{M/\max(\beta,1)} \int_{\beta y_{2}}^{M}
			\frac{a - y_{1}}{\left((a-y_{1})^{2} + y_{2}^{2}\right)^{1+\alpha/2}}
			\,dy_{1}\,dy_{2}
			= \int_{0}^{M/\max(\beta,1)}f(M,y_{2}) - f(\beta y_{2},y_{2})\,dy_{2}
		\end{align*}
		where we define
		\begin{align*}
			f(t,y_{2}) \coloneqq
			\begin{dcases}
				-\frac{1}{2}\ln\left((a-t)^{2} + y_{2}^{2}\right) & \textrm{if $\alpha=0$,} \\
				\frac{1}{\alpha}\frac{1}{\left((a-t)^{2} + y_{2}^{2}\right)^{\alpha/2}}
				& \textrm{otherwise.}
			\end{dcases}
		\end{align*}
		Differentiating twice by $a$, we obtain
		\begin{align*}
			A''(a) = \int_{0}^{M/\max(\beta,1)}
			\frac{(\alpha + 1)(a-M)^{2} - y_{2}^{2}}
			{\left((a-M)^{2} + y_{2}^{2}\right)^{2+\alpha/2}}
			+ \frac{y_{2}^{2} - (\alpha+1)(a - \beta y_{2})^{2}}
			{\left((a-\beta y_{2})^{2} + y_{2}^{2}\right)^{2+\alpha/2}}
			\,dy_{2} \eqqcolon A_{1}(a) + A_{2}(a).
		\end{align*}
		Since $a\leq\frac{M}{2}$, we have
		\[
			\abs{\frac{(\alpha + 1)(a-M)^{2} - y_{2}^{2}}
				{\left((a-M)^{2} + y_{2}^{2}\right)^{2+\alpha/2}}}
			\leq \frac{(\alpha+1)M^{2}}
			{(M/2)^{4+\alpha}},
		\]
		thus $\abs{A''(a) - A_{2}(a)} \leq C_\alp M^{-1-\alp}$.
		We further divide $A_{2}(a)$ into the sum of
		\begin{align*}
			A_{3}(a) &\coloneqq \int_{0}^{\infty}
			\frac{y_{2}^{2} - (\alpha+1)(a - \beta y_{2})^{2}}
			{\left((a-\beta y_{2})^{2} + y_{2}^{2}\right)^{2+\alpha/2}}\,dy_{2}
			\quad\textrm{and} \\
			A_{4}(a) &\coloneqq \int_{M/\max(\beta,1)}^{\infty}
			\frac{(\alpha+1)(a - \beta y_{2})^{2} - y_{2}^{2}}
			{\left((a-\beta y_{2})^{2} + y_{2}^{2}\right)^{2+\alpha/2}}\,dy_{2}.
		\end{align*}
		For the integrand in $A_{4}(a)$, if $\beta\geq 1$, then for $y_{2}$ in the domain of integration, we have $\beta y_{2} \geq M \geq 2a$,
		so $\beta y_{2} - a \geq \frac{\beta}{2}y_{2} \geq \frac{y_{2}}{2}$, which implies
		\[
			\abs{A_{4}(a)} \leq
			\int_{M/\beta}^{\infty}
			\frac{4}{(\beta y_{2} - a)^{2+\alpha}}\,dy_{2}
			\leq \frac{C_\alp}{M^{1+\alpha}}.
		\]
		On the other hand, if $\beta<1$, then for $y_{2}$ in the domain of integration, we have $\max(a,\beta y_{2})\leq y_{2}$, which implies
		\[
			\abs{A_{4}(a)} \leq
			\int_{M}^{\infty}
			\frac{\alpha+1}{y_{2}^{2+\alpha}}\,dy_{2}
			\leq \frac{C_\alp}{M^{1+\alpha}}.
		\]
		To evaluate $A_{3}(a)$, note that the denominator of the integrand of $A_{3}(a)$ is equal to
		\[
		\left( (\beta^{2}+1) \left(y_{2} - \frac{\beta a}{\beta^{2}+1}\right)^{2}
		+ \frac{a^{2}}{\beta^{2}+1}\right)^{2+\alpha/2},
		\]
		so by applying the change of variable $y_{2} = \frac{a(\beta + \tan\theta)}{\beta^{2}+1}$
		and the identity $a - \beta y_{2} = \frac{a(1 - \beta\tan\theta)}{\beta^{2}+1}$,
		we can rewrite $A_{3}(a)$ as
		\begin{align*}
			A_{3}(a) = \frac{1}{(\beta^{2}+1)^{1-\alpha/2}}\frac{1}{a^{1+\alpha}}
			\int_{-\tan^{-1}\beta}^{\pi/2}
			\left((\beta+\tan\theta)^{2} - (\alpha+1)(1-\beta\tan\theta)^{2}\right)
			\cos^{2+\alpha}\theta\,d\theta.
		\end{align*}
		The integrand can be rewritten as the product of $\cos^{\alpha}\theta$ and
		\begin{align*}
			&\left((\beta+\tan\theta)^{2} - (\alpha+1)(1-\beta\tan\theta)^{2}\right)\cos^{2}\theta \\
			&\quad\quad\quad
			= \left(1 - (\alpha+1)\beta^{2}\right)\sin^{2}\theta
			+ 2\beta(\alpha+2)\sin\theta\cos\theta
			+ (\beta^{2} - (\alpha+1))\cos^{2}\theta \\
			&\quad\quad\quad
			= (\alpha+2)\left(\frac{\beta^{2}-1}{2}\cos2\theta + \beta\sin2\theta\right)
			- \frac{\alpha(\beta^{2}+1)}{2}.
		\end{align*}
		By integration by parts, one can obtain the formulae
		\[
		(\alpha+2)\int_{s}^{t}\cos2\theta \cos^{\alpha}\theta\,d\theta
		- \alpha\int_{s}^{t}\cos^{\alpha}\theta\,d\theta
		= \left.2\sin\theta\cos^{\alpha+1}\theta\right|_{\theta=s}^{t}
		= \left.\frac{2\tan\theta}{(1+\tan^{2}\theta)^{1+\alpha/2}}\right|_{\theta=s}^{t}
		\]
		and  
		\[
		\int_{s}^{t}\sin2\theta \cos^{\alpha}\theta\,d\theta
		= 2\int_{\cos t}^{\cos s}u^{\alpha+1}\,du
		= \left.\frac{2u^{\alpha+2}}{\alpha+2}\right|_{u=\cos t}^{\cos s}
		= \left.\frac{2}{\alpha+2}\frac{1}{(1+\tan^{2}\theta)^{1+\alpha/2}}\right|_{\theta=t}^{s},
		\]
		so applying these we get
		\begin{align*}
			A_{3}(a) &= \frac{1}{(\beta^{2}+1)^{1-\alpha/2}}\frac{1}{a^{1+\alpha}}
			\left(
			\left.\frac{(\beta^{2}-1)u}{(1+u^{2})^{1+\alpha/2}}\right|_{u=-\beta}^{\infty}
			- \alpha\int_{-\tan^{-1}\beta}^{\pi/2}\cos^{\alpha}\theta\,d\theta
			+ \left.\frac{2\beta}{(1+u^{2})^{1+\alpha/2}}\right|_{u=-\beta}^{\infty}
			\right) \\
			&= \left(
			\frac{\beta}{\beta^{2}+1}
			- \frac{\alpha}{(\beta^{2}+1)^{1-\alpha/2}}
			\int_{-\tan^{-1}\beta}^{\pi/2}\cos^{\alpha}\theta\,d\theta
			\right)\frac{1}{a^{1+\alpha}}.
		\end{align*}
		Consequently, we arrive at the estimate \eqref{eq:vel-second-grad}. 
	\end{proof}

	\begin{proof}[Proof of Theorem~\ref{thm:main}]
		We assume the case $F(\alpha,\beta)>0$; the other case is almost identical. From \eqref{eq:kpp-evol}, we have \begin{equation*}
			\begin{split}
				\frac{\rd}{\rd t} \left( e^{2\int_0^t \rd_s v \cdot \bfT } \kpp \right) = -  e^{2\int_0^t \rd_s v \cdot \bfT } \rd_{s}^2 v \cdot \bfN . 
			\end{split}
		\end{equation*}
		Note that the estimate \eqref{eq:vel-second-grad} provides a uniform positive lower bound $b\leq \partial_{s}^{2}v_{2}((a,0),0)$ for all $(a,0)\in A$, where $A$ is defined as in \eqref{eq:annulus} with an appropriately modified $C_{\alpha}$. At $t=0$ and $x=(a,0)\in A$, we see that $\rd_{s}^2v\cdot \bfN= - \rd_{s}^2 v_{0,2} = -\rd_a^2 v_{0,2}\leq -b$. In particular, it follows that $\kpp(t,(a,0))$ is bounded below by a positive constant for all $(a,0)\in A$ and for all sufficiently small $t$, using the Taylor expansion of $\kpp$ in time. By symmetry, $\rd_{s}^2v\cdot \bfN$ must admit a uniform positive lower bound at $t=0$ on the intersection of $A$ and the upper edge, thus $\kpp$ is bounded above by a negative constant there for all sufficiently small $t$. Therefore, both of the edges must bent downward.
	\end{proof}
	
	\begin{proof}[Proof of Corollary~\ref{cor:main}]
		{We only need to ensure the existence of a single point where $\kpp$ becomes negative;} the idea is to simply regularize $\Omg_0$ in a way that its boundary is $C^\infty$-smooth everywhere and has strictly positive curvature except at a single point. To this end, we consider the following one-parameter family of modifications of $\Omg_0 $, denoted as $\Omg_0^{ \eps}$; we require that
		\begin{itemize}
			\item the curvature of $\partial\Omega_{0}^{\eps}$ is uniformly bounded in $\epsilon$, and is strictly positive except at a single fixed point $(a,0) \in \partial\Omg_0^{ \eps}\cap A$ where the curvature is zero,
			\item  $\partial\Omg_0^{ \eps}$ is tangent to  $\partial\Omg_0 $ at $(a,0)$, and
			\item $| \Omg_0^{ \eps} \triangle  \Omg_0 |<\eps$. 
		\end{itemize}
		It can be shown that, in the limit $\eps\to0$, \begin{equation*}
			\begin{split}
				\rd_{a}^2 v_{0,2}^{ \eps}(a,0) \to \rd_{a}^2 v_{0,2}(a,0).
			\end{split}
		\end{equation*}
		In particular, for all $\eps>0$ sufficiently small, we can arrange that \begin{equation*}
			\begin{split}
				\rd_{a}^2 v_{0,2}^{ \eps}(a,0) < 0,
			\end{split}
		\end{equation*}  where $v_{0,2}^{ \eps}$ is the vertical component of the velocity associated with the patch $\Omg_0^{ \eps}$. Then, we have that for all sufficiently small $t>0$, the curvature of $\partial\Omg^{\eps}(t)$ is negative somewhere.  This finishes the proof. 
	\end{proof}

	\begin{remark}\label{remk:velocity estimate}
		By performing a computation similar to what is shown in the proof of Proposition~\ref{prop:no strict convexity}, one can show that
		\[
			\abs{\partial_{s}v((a,0),0) + \frac{c_{\alpha}F(\alpha,\beta)}{\alpha a^{\alpha}}}
			\leq C_{\alpha}\left(\frac{1}{M^{\alpha}} + \frac{\delta^{2}}{a^{2+\alpha}}\right)
		\]
		holds for all $2\delta\leq a \leq M/2$ (with a possibly different $C_{\alpha}>0$). Therefore, when $0<\alpha<1$ and $\beta = \beta^{*}(\alpha)$ so that $F(\alpha,\beta)=0$, we obtain that $\partial_{s}v((\,\cdot\,,0),0) \to 0$ locally uniformly on $(0,\infty)$ when $M\to\infty$ and $\delta\to 0$. This means that when the corner becomes infinitely long and perfectly sharp, then the normal velocity on both of the edges of $\Omega_{0}$ must become a constant, which formally shows that the infinitely long and perfectly sharp corner with the angle $\beta^{*}(\alpha)$ is a steady solution to the gSQG equations with the parameter $\alpha$, relative to the frame moving together with the corner.
	\end{remark}
	
	\subsection*{Acknowledgment}
	
	\noindent JJ has been supported by NSF grant DMS-1900943. IJ has been supported by the Samsung Science and Technology Foundation under Project Number SSTF-BA2002-04.
	
	\bibliographystyle{amsplain}

\end{document}